\documentclass{gtart}

\def\ifplaintex{\expandafter\ifx\csname documentclass\endcsname\relax}

\def\gtp{{\mathsurround=0pt\it $\cal G\mskip-2mu$eometry \&\ 
$\cal T\!\!$opology $\cal P\!$ublications}}  

\def\Addressesr{\bigskip
{\small \parskip 0pt \leftskip 0pt \rightskip 0pt plus 1fil \def\\{\par}
\sl\theaddress\par
\medskip
\rm Email:\stdspace\tt\theemail\hfill\rm Received:\qua\receiveddate \par}}

\def\recd{{\small Received:\qua\receiveddate\ifx\reviseddate\relax
\else\qquad Revised:\qua\reviseddate\fi\par}} 

\def\lognumber#1{\def\thelognumber{#1}}
\def\volumenumber#1{\def\thevolumenumber{#1}}
\def\volumeyear#1{\def\thevolumeyear{#1}}
\def\papernumber#1{\def\thepapernumber{#1}}
\def\pagenumbers#1#2{\def\startpage{#1}\def\finishpage{#2}}
\def\published#1{\def\publishdate{#1}}

\def\received#1{\def\receiveddate{#1}}

\def\accepted#1{\def\accepteddate{#1}}

\def\asciiaddress#1{\def\theasciiaddress{#1}}

\long\def\asciiabstract#1{\long\def\theasciiabstract{#1}}

\let\\\par\let\thelognumber\relax\let\thevolumenumber\relax
\let\thepapernumber\relax\let\thevolumeyear\relax\let\startpage\relax
\let\finishpage\relax\let\publishdate\relax\let\receiveddate\relax
\let\reviseddate\relax\let\accepteddate\relax\let\theasciititle\relax
\let\theasciiauthors\relax\let\theasciiaddress\relax
\let\theasciiabstract\relax

\let\theasciiemail\relax

\ifplaintex
\font\logobig=cmssbx10 scaled 3836
\font\logomed=cmssbx10 scaled 2557
\else
\font\logobig=cmssbx10 scaled 4200
\font\logomed=cmssbx10 scaled 2800
\fi

\long\def\makeagttitle{   
\count0=\startpage
\agt\hfill      
\hbox to 45truept{\vbox to 0pt{\vglue -13truept{\logomed A\kern -.37em{\logobig 
T}\kern -.38em G}\vss}\hss}
\break
{\small Volume \thevolumenumber\ (\thevolumeyear)
\startpage--\finishpage\nl
Published: \publishdate}

\vglue .25truein

{\parskip=0pt\leftskip 0pt plus
1fil\def\\{\par\smallskip}{\Large\bf\thetitle}\par\medskip} \vglue
0.05truein

{\parskip=0pt\leftskip 0pt plus 1fil\def\\{\par}{\sc\theauthors}
\par\medskip}%
 
\vglue 0.03truein 

{\small\leftskip 25truept\rightskip 25truept{\bf Abstract}\stdspace\theabstract

{\bf AMS Classification}\stdspace\theprimaryclass
\ifx\thesecondaryclass\relax\else; \thesecondaryclass\fi\par
{\bf Keywords}\stdspace \thekeywords\par}\vglue 7truept

}   

\ifplaintex
\font\phead=cmsl9 scaled 950
\font\pnum=cmbx10 scaled 913
\font\pfoot=cmsl9 scaled 950
\headline{\vbox to 0pt{\vskip -4.5mm\line{\small\phead\ifnum
\count0=\startpage ISSN 1472-2739 (on-line) 1472-2747 (printed)
\hfill {\pnum\folio}\else\ifodd\count0\def\\{ }%
\ifx\theshorttitle\relax\thetitle\else\theshorttitle\fi\hfill{\pnum\folio}
\else\def\\{ and }{\pnum\folio}\hfill\ifx\theshortauthors\relax\theauthors
\else\theshortauthors\fi\fi\fi}\vss}}
\footline{\vbox to 0pt{\vglue 0mm\line{\small\pfoot\ifnum\count0=\startpage
\copyright\ \gtp\hfill\else
\agt, Volume \thevolumenumber\ (\thevolumeyear)\hfill\fi}\vss}}
\else
\font=cmsl9 scaled 1050
\font\lnum=cmbx10 
\font=cmsl9 scaled 1050
\makeatletter
\def\@oddhead{{\small\lhead\ifnum\count0=\startpage ISSN 1472-2739 
(on-line) 1472-2747 (printed)\hfill {\lnum\number\count0}\else\ifodd\count0
\def\\{ }\ifx\theshorttitle\relax \thetitle \else\theshorttitle\fi\hfill
{\lnum\number\count0}\else\def\\{ and }{\lnum\number\count0}
\hfill\ifx\theshortauthors\relax 
\theauthors\else\theshortauthors\fi\fi\fi}}\def\@evenhead{\@oddhead}
\def\@oddfoot{\small\lfoot\ifnum\count0=\startpage\copyright\ \gtp\hfill\else
\agt, Volume \thevolumenumber\ (\thevolumeyear)\hfill\fi}
\def\@evenfoot{\@oddfoot}
\makeatother
\fi
\let\maketitlepage\makeagttitle
\let\makeshorttitle\maketitlepage
\let\maketitle\maketitlepage

\newwrite\gtoutfile
\long\gdef\makeheadfile{  
{\def\\{, }\def\s{ }
\immediate\openout\gtoutfile head.xxx
\immediate\write\gtoutfile{To: math@arxiv.org}
\immediate\write\gtoutfile{Subject: put OR rep NNNNN:ppppp}
\immediate\write\gtoutfile{--text follows this line--}
\immediate\write\gtoutfile{Proxy-for: \ifx\theasciiauthors\relax
\theauthors\else\theasciiauthors\fi\s<\ifx\theasciiemail\relax\theemail\else\theasciiemail\fi>}
\immediate\write\gtoutfile{\noexpand\\}
\immediate\write\gtoutfile{Authors: \ifx\theasciiauthors\relax
\theauthors\else\theasciiauthors\fi}
{\def\\{ }\immediate\write\gtoutfile{Title: \ifx\theasciititle\relax
\thetitle\else\theasciititle\fi}}
\immediate\write\gtoutfile{Subj-class: GT or SG, GR etc}
\immediate\write\gtoutfile{MSC-class: \theprimaryclass\ifx\thesecondaryclass\relax\else, \thesecondaryclass\fi}
\immediate\write\gtoutfile{Journal-ref: Algebr. Geom. Topol. \thevolumenumber\s
(\thevolumeyear) \startpage-\finishpage}
\immediate\write\gtoutfile{Comments: Published by Algebraic and
Geometric Topology at}
\immediate\write\gtoutfile{\s\s\s  http://www.maths.warwick.ac.uk/agt/AGTVol\thevolumenumber/agt-\thevolumenumber-\thepapernumber.abs.html}
\immediate\write\gtoutfile{\noexpand\\}
\immediate\write\gtoutfile{}
\ifx\theasciiabstract\relax
\immediate\write\gtoutfile{\theabstract}\else
\immediate\write\gtoutfile{\theasciiabstract}\fi
\immediate\write\gtoutfile{}
\immediate\write\gtoutfile{\noexpand\\}
\immediate\write\gtoutfile{}
\immediate\closeout\gtoutfile}}  

\def\maketitlepage{\makeagttitle\makeheadfile}
\let\makeshorttitle\maketitlepage
\let\maketitle\maketitlepage

\lognumber{17}
\volumenumber{3}
\volumeyear{2003}
\papernumber{17}
\published{16 June 2003}
\pagenumbers{537}{556}
\received{12 December 2002}
\accepted{22 May 2003}

\usepackage{amsmath,latexsym,amssymb,pstricks,mathrsfs}

\newtheorem{thm}{Theorem}[section]
\newtheorem{lem}[thm]{Lemma}
\newtheorem{prop}[thm]{Proposition}

\newtheorem{rem}[thm]{Remark}

\newcommand{\BZ}{{\mathbb{Z}}}

\newcommand{\BR}{{\mathscr{R}}}
\newcommand{\BS}{{\mathscr{S}}}

\newcommand{\BB}{{\mathcal{B}}}
\newcommand{\BK}{{\mathcal{K}}}

\newcommand{\bo}{{\mathfrak o}}
\newcommand{\fsl}{{\mathfrak s \mathfrak l}}

\def\mapright#1{\smash{\mathop{\longrightarrow}\limits^{#1}}}

\newcommand\mychoose[2]{\genfrac{(}{)}{0pt}{}{#1}{#2}}
\newcommand\qchoose[2]{\genfrac{[}{]}{0pt}{}{#1}{#2}}

\newcommand{\BOX}{
\psset{unit=.35cm}
\begin{pspicture}[.4](-.3,-2)(5,3)
\small
\psframe(0,0)(4.5,1)
\psline(1,0)(1,-2)
\psline(1.5,0)(1.5,-2)
\psline(3.5,0)(3.5,-2)
\psline(1,1)(1,3)
\psline(3.5,1)(3.5,3)
\psline(1.5,1)(1.5,3)
\put(2,1.8){$\cdots$}
\put(2,-.8){$\cdots$}
\end{pspicture}
}

\newcommand{\LINE}{
\psset{unit=.35cm}
\begin{pspicture}[.4](-.3,-2)(1,3)
\small
\psline(0,-2)(0,3)
\put(.3,1.8){$n$}
\end{pspicture}
}

\newcommand{\VERTEX}{
\psset{unit=.5cm}
\begin{pspicture}[.4](-2,-2)(2,1.5)
\small
\psline(0,0)(0,-2)
\psline(0,0)(-2,1.4)
\psline(0,0)(2,1.4)
\put(-1.6,.4){$a$}
\put(1.4,.4){$b$}
\put(.2,-1.4){$c$}
\end{pspicture}
}

\newcommand{\VERTEXi}{
\psset{unit=.5cm}
\begin{pspicture}[.4](-2,-2)(2,1.5)
\small
\psline(0,-1)(0,-2)
\psline(-1.2,.8)(-2,1.4)
\psline(1.2,.8)(2,1.4)
\psframe(-.5,-1)(.5,-.5)
\rput{120}(1,.7){\psframe(-.5,-.25)(.5,.25)}
\rput{240}(-1,.7){\psframe(-.5,-.25)(.5,.25)}
\psbezier[showpoints=false](-.2,-.5)(-.3,0)(0,0)(-.9,.4)
\psbezier[showpoints=false](.2,-.5)(.3,0)(0,0)(.9,.4)
\pscurve(-.72,.72)(0,.4)(.72,.72)
\put(-2,.7){$a$}
\put(1.8,.7){$b$}
\put(.2,-1.7){$c$}
\put(-.1,.7){$k$}
\put(.6,-.3){$i$}
\put(-.9,-.3){$j$}
\end{pspicture}
}

\newcommand{\FUS}{
\psset{unit=.5cm}
\begin{pspicture}[.25](0,0)(3,2)
\small
\put(.5,.2){$b$}
\put(.5,1.2){$a$}
\psline(0,0)(3,0)
\psline(0,1)(3,1)
\end{pspicture}
}

\newcommand{\FUSi}{
\psset{unit=.5cm}
\begin{pspicture}[.25](0,0)(4,2)
\small
\put(.5,-.3){$b$}
\put(.4,1){$a$}
\put(1.8,0){$c$}
\put(3.35,-.4){$b$}
\put(3.4,1){$a$}
\psline(0,0)(1,.5)
\psline(3,.5)(4,0)
\psline(1,.5)(3,.5)
\psline(0,1)(1,.5)
\psline(4,1)(3,.5)
\end{pspicture}
}

\newcommand{\FUSa}{
\psset{unit=.5cm}
\begin{pspicture}[.25](0,0)(3,2)
\small
\put(.5,.2){$n$}
\put(.5,1.2){$n$}
\psline(0,0)(3,0)
\psline(0,1)(3,1)
\end{pspicture}
}

\newcommand{\FUSai}{
\psset{unit=.5cm}
\begin{pspicture}[.25](0,0)(4,2)
\small
\put(.5,-.3){$n$}
\put(.4,1){$n$}
\put(1.7,-.1){$2n$}
\put(3.4,-.3){$n$}
\put(3.35,1.05){$n$}
\psline(0,0)(1,.5)
\psline(3,.5)(4,0)
\psline(1,.5)(3,.5)
\psline(0,1)(1,.5)
\psline(4,1)(3,.5)
\end{pspicture}
}

\newcommand{\DEL}{
\psset{unit=.5cm}
\begin{pspicture}[.4](0,0)(3,1)
\small
\psline(0,0)(.3,.3)
\psline(.7,.7)(1,1)(1.5,.5)
\psline(0,1)(1,0)(1.5,.5)
\psline(1.5,.5)(2.5,.5)
\put(.2,-.4){$a$}
\put(.2,1.1){$b$}
\put(1.8,0){$c$}
\end{pspicture}
}
\newcommand{\DELi}{
\psset{unit=.5cm}
\begin{pspicture}[.4](0,-.5)(3,1.5)
\small
\psline(0,0)(1.5,.5)
\psline(0,1)(1.5,.5)
\psline(1.5,.5)(2.5,.5)
\put(.2,-.4){$a$}
\put(.2,1.1){$b$}
\put(1.8,0){$c$}
\end{pspicture}
}

\newcommand{\THETA}{
\psset{unit=.5cm}
\begin{pspicture}[.3](-.5,-1)(3.5,1.5)
\small
\psline(0,0)(3,0)
\psellipse(1.5,0)(1.5,1)
\put(1.2,-.8){$a$}
\put(1.2,1.2){$c$}
\put(1.2,.2){$b$}
\end{pspicture}
}

\newcommand{\CAL}{
\psset{unit=.5cm}
\begin{pspicture}[.45](-1,0)(2,5)
\small
\psline(0,0)(0,1)(1,2.5)(0,4)(0,5)
\psline(0,1)(-1,2.5)(0,4)
\psline(-1,2.5)(1,2.5)
\put(.2,.1){$2n$}
\put(.2,4.4){$2n$}
\put(.7,1.3){$n$}
\put(.7,3.2){$n$}
\put(-1,1.3){$n$}
\put(-1.1,3.2){$n$}
\put(-.4,2.7){$2k$}
\end{pspicture}
}

\newcommand{\CALa}{
\psset{unit=.5cm}
\begin{pspicture}[.45](-1,0)(2,5)
\small
\psline(0,0)(0,1)(1,2.5)(0,4)(0,5)
\psline(0,1)(-1,2.5)(0,4)
\psline(-1,2.5)(1,2.5)
\put(.2,.1){$2n$}
\put(.2,4.4){$2n$}
\put(.7,1.3){$n$}
\put(.7,3.2){$n$}
\put(-1,1.3){$n$}
\put(-1.1,3.2){$n$}
\put(-.4,2.7){$2n$}
\end{pspicture}
}

\newcommand{\CALip}{
\psset{unit=.5cm}
\begin{pspicture}[.45](-1,0)(2,5)
\small
\psline(0,0)(0,1)(1,2.5)(0,4)(0,5)
\psline(0,1)(-1,2.5)(0,4)
\put(.2,.1){$2n$}
\put(.2,4.4){$2n$}
\put(.7,3.2){$n$}
\put(-1.1,3.2){$n$}
\end{pspicture}
}

\newcommand{\CALo}{
\psset{unit=.5cm}
\begin{pspicture}[.45](-1,0)(1,5)
\small
\psline(0,0)(0,5)
\put(.2,4.4){$2n$}
\end{pspicture}
}

\newcommand{\CALiio}{
\psline(0,0)(0,1)
\psline(0,4)(0,5)
\put(.2,.1){$2n$}
\put(.2,4.4){$2n$}
\psarc(0,2){1}{180}{0}
\psarc(0,3){1}{0}{180}
\psline(-1,3)(-.5,3)\psline(-.1,3)(1,3)
\psline(-1,2)(-.5,2)\psline(-.1,2)(1,2)
\psarc(0,1.8){.3}{180}{0}
\psarc(0,3.2){.3}{0}{180}
\psline(-.3,1.8)(-.3,3.2)
\psline(.3,2.2)(.3,2.8)
\put(-1.6,3.3){$n$}
\put(1.1,1.5){$n$}
}
\newcommand{\CALii}{
\psset{unit=.5cm}
\begin{pspicture}[.45](-2,0)(2,5)
\small
\CALiio
\put(-1.4,2.4){$R_n$}
\end{pspicture}
}

\newcommand{\CALiia}{
\psset{unit=.5cm}
\begin{pspicture}[.45](-2,0)(2,5)
\CALiio
\put(-1.4,2.4){$R_k$}
\end{pspicture}
}

\newcommand{\CALiic}{
\psset{unit=.5cm}
\begin{pspicture}[.45](-2,0)(2,5)
\small
\CALiio
\put(-1.3,2.3){$\omega^p$}
\end{pspicture}
}

\newcommand{\CALiii}{
\psset{unit=.5cm}
\begin{pspicture}[.45](-3,0)(3,5)
\small
\psline(0,0)(0,1)
\psline(0,4)(0,5)
\put(.2,.1){$2n$}
\put(.2,4.4){$2n$}
\psarc(-1,2){1}{180}{270}\psarc(1,2){1}{270}{0}
\psarc(-1,3){1}{90}{180}\psarc(1,3){1}{0}{90}
\psline(-1,4)(1,4)
\psline(-1,1)(1,1)
\put(-2.6,3.2){$n$}
\put(2.1,1.5){$n$}
\psline(-2,3)(-2,2.5)(-.5,2.5)\psline(-.1,2.5)(2,2.5)(2,3)
\psline(-2,2)(-2,2.5)
\psline(2,2.5)(2,2)
\psarc(0,1.8){.3}{180}{0}
\psarc(0,3.2){.3}{0}{180}
\psline(-.3,1.8)(-.3,3.2)
\psline(.3,1.8)(.3,2.3)\psline(.3,2.7)(.3,3.2)
\put(-1.4,2.9){$R_n$}
\put(.8,1.9){$2n$}
\end{pspicture}
}

\newcommand{\CALiv}{
\psset{unit=.5cm}
\begin{pspicture}[.45](-3,0)(3,5)
\small
\psline(0,0)(0,1)
\psline(0,4)(0,5)
\put(.2,.1){$2n$}
\put(.2,4.4){$2n$}
\psarc(-1,2){1}{180}{270}\psarc(1,2){1}{270}{0}
\psarc(-1,3){1}{90}{180}\psarc(1,3){1}{0}{90}
\psline(-1,4)(1,4)
\psline(-1,1)(1,1)
\put(-2.6,3.2){$n$}
\put(2.1,1.5){$n$}
\psline(-2,3)(-1,2)
\psline(-2,2)(-1.7,2.3)
\psline(-1.3,2.7)(-1,3)
\psline(1,3)(2,2)
\psline(1,2)(1.3,2.3)
\psline(1.7,2.7)(2,3)
\put(-.5,2.4){$\cdots$}
\end{pspicture}
}

\newcommand{\TWi}{
\psset{unit=.5cm}
\begin{pspicture}[.45](-5,0)(2,5)
\small
\psarc(0,2){1}{180}{0}
\psarc(0,3){1}{0}{180}
\psline(-1,3)(0,2)(.3,2.3)
\psline(.7,2.7)(1,3)
\psline(1,2)(0,3)(-.3,2.7)
\psline(-.7,2.3)(-1,2)
\put(1.1,3.2){$n$}
\put(1.1,1.5){$n$}
\psline(-3,1)(-3,1.8)\psline(-3,2.2)(-3,4)
\scalebox{1 .5}{\psarc(-3,5){1}{110}{70}}
\psecurve(-3,3.8)(-3,4)(-3,4.2)(-2,4.8)(-1,4.8)(0,4.2)(0,4)(0,3.8)
\psecurve(-3,1.2)(-3,1)(-3,.8)(-2,.2)(-1,.2)(0,.8)(0,1)(0,1.2)
\put(-1.5,4.2){$2n$}
\put(-2.5,1.4){$R'_n$}
\end{pspicture}
}

\newcommand{\TW}{
\psset{unit=.5cm}
\begin{pspicture}[.45](-5,-.5)(2,5)
\small
\rput(0,1){
\psline(-1,3)(0,2)(.3,2.3)
\psline(.7,2.7)(1,3)
\psline(1,2)(0,3)(-.3,2.7)
\psline(-.7,2.3)(-1,2)
}
\rput[bl]{90}(-1,3){
\psline(-1,3)(0,2)(.3,2.3)
\psline(.7,2.7)(1,3)
\psline(1,2)(0,3)(-.3,2.7)
\psline(-.7,2.3)(-1,2)
}
\put(-3.6,1.2){$\cdot$}\put(-3.6,1.5){$\cdot$}\put(-3.6,1.8){$\cdot$}
\pscurve(-1,4)(-2,4.5)(-3,4)
\pscurve(1,4)(.9,4.4)(-2,5)(-3.9,4.4)(-4,4)
\pscurve(-1,1)(-2,0.5)(-3,1)
\pscurve(1,1)(.9,.6)(-2,0)(-3.9,.6)(-4,1)
\psline(-1,1)(-1,3)
\psline(1,1)(1,3)
\put(-6.2,2.8){$p$ full}
\put(-6.2,2){twists}
\end{pspicture}
}

\newcommand{\TWao}{
\rput(0,1){
\psline(-1,3)(0,2)(.3,2.3)
\psline(.7,2.7)(1,3)
\psline(1,2)(0,3)(-.3,2.7)
\psline(-.7,2.3)(-1,2)
}
\pscurve(-1,4)(-2,4.5)(-3,4)
\pscurve(1,4)(.9,4.4)(-2,5)(-3.9,4.4)(-4,4)
\pscurve(-1,1)(-2,0.5)(-3,1)
\pscurve(1,1)(.9,.6)(-2,0)(-3.9,.6)(-4,1)
\psline(-1,1)(-1,3)
\pscurve(1,1)(1.2,1.6)(1.5,1.7)(1.7,1.5)(1.5,1.3)(1.3,1.4)
\psline(1.1,1.8)(1,2)
\rput(0,1){
\pscurve(1,1)(1.2,1.6)(1.5,1.7)(1.7,1.5)(1.5,1.3)(1.3,1.4)
\psline(1.1,1.8)(1,2)
}
\psline(-4,4)(-4,2.2)\psline(-4,1.7)(-4,1)
\psline(-3,4)(-3,2.2)\psline(-3,1.7)(-3,1)
\scalebox{1 .5}{\psarc(-3.5,5){1.2}{130}{50}}
\scalebox{1 .5}{\psarc(-3.5,5){1.2}{85}{110}}
}
\newcommand{\TWa}{
\psset{unit=.5cm}
\begin{pspicture}[.45](-7,-.5)(4,5)
\small
\TWao
\put(-6,2.3){$\omega^p$}
\end{pspicture}
}

\newcommand{\TWan}{
\psset{unit=.5cm}
\begin{pspicture}[.45](-7,-.5)(2,5)
\small
\TWao
\put(1,4){$R_n$}
\put(-6.1,2.3){$R'_k$}
\end{pspicture}
}

\newcommand{\TWnn}{
\psset{unit=.5cm}
\begin{pspicture}[.45](-6.5,0)(2,5)
\small
\TWao
\put(1,4){$R_n$}
\put(-6.1,2.3){$R'_n$}
\end{pspicture}
}

\newcommand{\TWnnn}{
\psset{unit=.5cm}
\begin{pspicture}[.45](-7,-.5)(2,5)
\TWao
\small
\put(1,4){$e_n$}
\put(-6.1,2.3){$R'_n$}
\end{pspicture}
}

\newcommand{\BOXES}{
\psframe(-.4,-.4)(.4,0)
\psframe(-.4,2)(.4,2.4)
\psframe(3.6,-.4)(4.4,0)
\psframe(3.6,2)(4.4,2.4)
}

\newcommand{\REC}{
\psset{unit=.5cm}
\begin{pspicture}[.45](-1,-1)(5,3)
\small
\BOXES
\psline(0,2)(2,0)(2.7,.7)
\psline(3.3,1.3)(4,2)
\psline(0,0)(.7,.7)
\psline(1.3,1.3)(2,2)(4,0)
\put(-.4,1.5){$n$}
\put(-.4,.3){$n$}
\end{pspicture}
}

\newcommand{\RECp}{
\psset{unit=.5cm}
\begin{pspicture}[.45](-1,-1)(7,3)
\small
\psframe(-.4,-.4)(.4,0)
\psframe(-.4,2)(.4,2.4)
\psframe(5.6,-.4)(6.4,0)
\psframe(5.6,2)(6.4,2.4)
\psline(0,2)(2,0)\psline(4,0)(4.7,.7)
\psline(5.3,1.3)(6,2)
\psline(0,0)(.7,.7)
\psline(1.3,1.3)(2,2)\psline (4,2)(6,0)
\put(-.4,1.5){$n$}
\put(-.4,.3){$n$}
\put(2.6,.8){$\cdots$}
\end{pspicture}
}

\newcommand{\RECj}{
\psset{unit=.5cm}
\begin{pspicture}[.45](-1,-1)(5,3)
\small
\BOXES
\psline(0,2)(2,0)(2.7,.7)
\psline(3.3,1.3)(4,2)
\psline(0,0)(.7,.7)
\psline(1.3,1.3)(2,2)(4,0)
\put(-.4,1.5){$n$}
\put(-.4,.3){$p$}
\end{pspicture}
}

\newcommand{\RECi}{
\psset{unit=.5cm}
\begin{pspicture}[.45](-1,-1)(5,3)
\small
\psframe(-.4,2)(.4,2.4)
\psframe(3.6,2)(4.4,2.4)
\psline(0,2)(2,0)(2.7,.7)
\psline(3.3,1.3)(4,2)
\psline(0,0)(.7,.7)
\psline(1.3,1.3)(2,2)(4,0)
\put(-.4,1.5){$n$}
\end{pspicture}
}

\newcommand{\RECa}{
\psset{unit=.5cm}
\begin{pspicture}[.45](-1,-1)(5,3)
\small
\BOXES
\pscurve(.2,0)(.3,.1)(2,.4)(3.7,.1)(3.8,0)
\psline(-.2,0)(-.2,2)
\pscurve(.2,2)(.3,1.9)(2,1.6)(3.7,1.9)(3.8,2)
\psline(4.2,0)(4.2,2)
\put(0,.8){$k$}
\put(3.6,.8){$k$}
\put(1.2,1.9){$n-k$}
\put(1.2,-.4){$n-k$}
\end{pspicture}
}

\newcommand{\RECaj}{
\psset{unit=.5cm}
\begin{pspicture}[.45](-1,-1)(5,3)
\small
\BOXES
\pscurve(.2,0)(.3,.1)(2,.4)(3.7,.1)(3.8,0)
\psline(-.2,0)(-.2,2)
\pscurve(.2,2)(.3,1.9)(2,1.6)(3.7,1.9)(3.8,2)
\psline(4.2,0)(4.2,2)
\put(0,.8){$k$}
\put(3.6,.8){$k$}
\put(1.2,1.9){$n-k$}
\put(1.2,-.3){$p-k$}
\end{pspicture}
}

\newcommand{\RECii}{
\psset{unit=.5cm}
\begin{pspicture}[.45](-1,-1)(5,3)
\small
\psframe(-.4,2)(.4,2.4)
\psframe(3.6,2)(4.4,2.4)
\pscurve(.2,0)(.3,.1)(2,.4)(3.7,.1)(3.8,0)
\pscurve(0,2)(.1,1.9)(2,1.6)(3.9,1.9)(4,2)
\put(1.5,1.9){$n$}
\end{pspicture}
}

\newcommand{\RECiii}{
\psset{unit=.5cm}
\begin{pspicture}[.45](-1,-1)(5,3)
\small
\psframe(-.4,2)(.4,2.4)
\psframe(3.6,2)(4.4,2.4)
\psline(-.2,0)(-.2,2)
\pscurve(.2,2)(.3,1.9)(2,1.6)(3.7,1.9)(3.8,2)
\psline(4.2,0)(4.2,2)
\put(1.2,1.9){$n-1$}
\end{pspicture}
}

\newcommand{\RECb}{
\psset{unit=.5cm}
\begin{pspicture}[.45](-1,-1)(9,3)
\small
\psframe(-.4,-.4)(.4,0)
\psframe(-.4,2)(.4,2.4)
\psframe(4,-.4)(4.4,.4)
\psframe(4,1.6)(4.4,2.4)
\psframe(7.6,-.4)(8.4,0)
\psframe(7.6,2)(8.4,2.4)
\psline(0,2)(2,0)(2.7,.7)
\psline(3.3,1.3)(4,2)
\psline(0,0)(.7,.7)
\psline(1.3,1.3)(2,2)(4,0)
\psarc(4.4,1){.8}{270}{90}
\pscurve(4.4,2.2)(4.5,2.2)(6,1.6)(7.7,1.9)(7.8,2)
\pscurve(4.4,-.2)(4.5,-.2)(6,.4)(7.7,.1)(7.8,0)
\psline(8.2,0)(8.2,2)
\put(-.4,1.5){$n$}
\put(-.4,.3){$n$}
\put(5.3,.8){$k$}
\put(8.4,.8){$k$}
\put(5.5,1.9){$n-k$}
\put(5.5,-.4){$n-k$}
\end{pspicture}
}

\newcommand{\RECc}{
\psset{unit=.5cm}
\begin{pspicture}[.45](-1,-1)(5,3)
\small
\BOXES
\psline(-.2,0)(-.2,2)
\psline(4.2,0)(4.2,2)
\psline(.2,2)(2,.2)(2.6,.8)
\psline(3.1,1.3)(3.8,2)
\psline(.2,0)(.9,.7)
\psline(1.4,1.2)(2,1.8)(3.8,0)
\put(0,.8){$k$}
\put(4.4,.8){$k$}
\put(1.2,2.1){$n-k$}
\put(1.2,-.4){$n-k$}
\end{pspicture}
}

\newcommand{\RECd}{
\psset{unit=.5cm}
\begin{pspicture}[.45](-1,-1)(9,3)
\small
\psframe(-.4,-.4)(.4,0)
\psframe(-.4,2)(.4,2.4)
\psframe(7.6,-.4)(8.4,0)
\psframe(7.6,2)(8.4,2.4)

\psline(.2,2)(2,.2)(2.6,.8)
\psline(3.1,1.3)(4,2.2)
\psline(.2,0)(.9,.7)
\psline(1.4,1.2)(2,1.8)(4,-.2)

\psline(-.2,2)(2,-.2)(2.8,.6)
\psline(3.4,1.2)(4,1.8)
\psline(-.2,0)(.6,.8)
\psline(1.2,1.4)(2,2.2)(4,.2)

\psarc(4,1){.8}{270}{90}
\pscurve(4,2.2)(4.4,2.2)(4.5,2.2)(6,1.6)(7.7,1.9)(7.8,2)
\pscurve(4,-.2)(4.4,-.2)(4.5,-.2)(6,.4)(7.7,.1)(7.8,0)
\psline(8.2,0)(8.2,2)
\put(5,.8){$k$}
\put(8.4,.8){$k$}
\put(5.5,1.9){$n-k$}
\put(5.5,-.4){$n-k$}
\end{pspicture}
}

\newcommand{\RECe}{
\psset{unit=.5cm}
\begin{pspicture}[.45](-1,-1)(5,3)
\small
\BOXES
\psline(4.2,0)(4.2,2)
\psline(.2,2)(.35,1.85)\psline(.6,1.6)(2,.2)(2.6,.8)
\psline(3.1,1.3)(3.8,2)
\psline(.2,0)(.6,.4)\psline(.85,.65)(1.05,.85)
\psline(1.4,1.2)(2,1.8)(3.8,0)
\pscurve(-.2,0)(0,.3)(.3,.3)
\pscurve(.7,.2)(.9,.2)(.9,.4)(.4,.6)(.2,1)(.4,1.4)(.5,1.4)
\pscurve(.9,1.6)(.9,1.8)(.7,1.8)(.3,1.7)(0,1.7)(-.2,2)
\put(-.4,.8){$k$}
\put(4.4,.8){$k$}
\put(1.2,2.1){$n-k$}
\put(1.2,-.4){$n-k$}
\end{pspicture}
}

\newcommand{\KAUF}{
\psset{unit=.5cm}
\begin{pspicture}[.45](0,0)(2,2)
\small
\psline(0,2)(2,0)
\psline(0,0)(.7,.7)
\psline(1.3,1.3)(2,2)
\end{pspicture}
}

\newcommand{\KAUFi}{
\psset{unit=.5cm}
\begin{pspicture}[.45](0,0)(2,2)
\scalebox{1 .5}{\psarc(1,0){1}{0}{180}}
\scalebox{1 .5}{\psarc(.8,4){1}{180}{0}}
\end{pspicture}
}

\newcommand{\KAUFii}{
\psset{unit=.5cm}
\begin{pspicture}[.45](0,0)(2,2)
\scalebox{.5 1}{\psarc(0,1){1}{270}{90}}
\scalebox{.5 1}{\psarc(3.5,1){1}{90}{270}}
\end{pspicture}
}

\newcommand{\KAUFiii}{
\psset{unit=.5cm}
\begin{pspicture}[.45](0,0)(2,2)
\pscircle(1,1){.8}
\end{pspicture}
}

\newcommand{\IDEM}{
\psset{unit=.4cm}
\begin{pspicture}[.45](-1,-.5)(5.5,2.5)
\small
\psframe(1.6,0)(2.2,2)
\psframe(3.2,1)(3.6,2)
\psline(-1,1)(1.6,1)
\psline(2.2,.5)(4.5,.5)
\psline(2.2,1.5)(3.2,1.5)\psline(3.6,1.5)(4.5,1.5)
\put(-.8,1.4){$p+q$}
\put(4.8,1.5){$p$}
\put(4.8,.2){$q$}
\end{pspicture}
}

\newcommand{\IDEMii}{
\psset{unit=.4cm}
\begin{pspicture}[.45](-1,-.5)(5.5,2.5)
\small
\psframe(1.6,0)(2.2,2)
\psline(-1,1)(1.6,1)
\psline(2.2,.5)(4.5,1.5)
\psline(2.2,1.5)(3,1.15)\psline(3.7,.8)(4.5,.5)
\put(-.8,1.4){$p+q$}
\put(4.8,1.5){$p$}
\put(4.8,.2){$q$}
\end{pspicture}
}

\newcommand{\IDEMiii}{
\psset{unit=.4cm}
\begin{pspicture}[.45](-1,-.5)(5.5,2.5)
\small
\psframe(1.6,0)(2.2,2)
\psline(-1,1)(1.6,1)
\psline(2.2,.5)(4.5,.5)
\psline(2.2,1.5)(4.5,1.5)
\put(-.8,1.4){$p+q$}
\put(4.8,1.5){$p$}
\put(4.8,.2){$q$}
\end{pspicture}
}

\newcommand{\CIRC}{
\psset{unit=.5cm}
\begin{pspicture}[.45](-1,0)(3,3)
\small
\psline(1,0)(1,.8)\psline(1,1.4)(1,3)
\scalebox{1 .5}{\psarc(1,3){.8}{110}{70}}
\put(1.3,2.6){$x_{even}$}
\put(-1,1.5){$\omega_+$}
\end{pspicture}
}

\newcommand{\CIRCi}{
\psset{unit=.5cm}
\begin{pspicture}[.45](0,0)(2,3)
\small
\put(1.3,2.6){$x_{even}$}
\pscurve(1,1)(1.2,1.6)(1.5,1.7)(1.7,1.5)(1.5,1.3)(1.3,1.4)
\psline(1.1,1.8)(1,2)(1,3)
\psline(1,0)(1,1)
\end{pspicture}
}

\begin{document}

\title{Skein-theoretical derivation\\of some formulas of Habiro}

\author{ Gregor Masbaum}
\address{Institut de Math{\'e}matiques de Jussieu (UMR 7586 du CNRS)\\
Universit{\'e} Paris 7 (Denis Diderot), Case 7012\\
2, place Jussieu, 75251 Paris Cedex 05, France}
\asciiaddress{Institut de Mathematiques de Jussieu (UMR 7586 du CNRS)\\
Universite Paris 7 (Denis Diderot), Case 7012\\2, place Jussieu, 75251 
Paris Cedex 05, France}
\email{masbaum@math.jussieu.fr}

\begin{abstract}  We use skein theory  to compute the coefficients of certain power series considered by Habiro in his theory of $\fsl_2$ invariants of integral homology $3$-spheres. Habiro originally derived these formulas using the quantum group $U_q\fsl_2$. As an application, we give a formula for the colored Jones polynomial of twist knots, generalizing formulas of Habiro and Le for the trefoil and the figure eight knot.
\end{abstract}
\asciiabstract{We use skein theory to compute the coefficients of
certain power series considered by Habiro in his theory of sl_2
invariants of integral homology 3-spheres. Habiro originally derived
these formulas using the quantum group U_q sl_2. As an application,
we give a formula for the colored Jones polynomial of twist knots,
generalizing formulas of Habiro and Le for the trefoil and the figure
eight knot.}


\maketitle

\section*{Introduction}

In a talk at the Mittag-Leffler Institute in May 1999, K.~Habiro announced 
a new approach to computing the colored Jones polynomial of knots and quantum
$\fsl_2$ invariants of integral homology $3$-spheres. For an exposition, see 
his paper \cite{Hab2} (some results are already announced in \cite{Hab}). His invariant for homology spheres  recovers both the $sl_2$ Reshetikhin-Turaev invariants at roots of unity, and Ohtsuki's power series invariants.  
Later, Habiro and T.Q.T.~Le generalized this to all quantum invariants associated to simple Lie algebras.

In the $\fsl_2$ case, quantum invariants can be expressed in terms of skein theory, using the Jones polynomial or the Kauffman bracket. Habiro's invariant
for homology spheres can be constructed using certain skein elements $\omega=\omega_+$ and $\omega^{-1}=\omega_-$ such that circling an even number of strands with $\omega_+$ (resp. $\omega_-$) induces a positive (resp. negative) full twist: 
\begin{equation}\label{fig-1}
\CIRC = \CIRCi
\end{equation}
More precisely, $\omega_+$ and $\omega_-$ are not skein elements, but infinite sums ({\em i.e.} power series) of such. But as long as they encircle an even number of strands (corresponding to a strand colored by an integer-spin representation of $\fsl_2$), the result is well-defined. Also, it makes sense to consider  powers of $\omega=\omega_+$, and circling an even number of strands with  $\omega^p$  induces $p$ positive  full twists, where $p\in \BZ$. 

The main purpose of this paper is to give skein-theoretical proofs of Habiro's
formulas for $\omega_+$ and $\omega_-$ (they are stated already in \cite{Hab}) and for $\omega^p$ (this formula  will  appear in \cite{Hab3}). Habiro's original proofs of these formulas use the quantum group $U_q\fsl_2$.

This paper is organized as follows. After stating Habiro's formula for $\omega_+$ and $\omega_-$ in Section~1, we give a proof using orthogonal polynomials along the lines of \cite{BHMV1} in Section~2. This proof is quite straightforward, although the computations are a little bit more involved than in \cite{BHMV1}. Unfortunately, it seems difficult to use this approach to compute the coefficients of $\omega^p$ for $|p|\geq 2$. Therefore, in Section~3 we start afresh using the Kauffman bracket graphical calculus. A first expression for $\omega^p$ in Theorem~\ref{th2} is easily obtained, but it is not quite good enough, as 
an important divisibility property of the coefficients of $\omega^p$ is not clear from this formula.  This property is then shown in Section~4 by some more skein theory. The final expression for $\omega^p$ obtained in Theorem~\ref{th3} is equivalent to Habiro's one   from \cite{Hab3}. (The results of Section~2 are not used here, so that this gives an independent proof in the $p=\pm 1$ case as well.) 

To illustrate one use of $\omega^p$, we conclude the paper in Section~5 by giving a formula for the colored Jones polynomial of twist knots. This generalizes formulas of Habiro \cite{Hab} (see also Le \cite{Le1,Le2}) for the trefoil and the figure eight knot. (For those two knots, one only needs Habiro's original $\omega_+$ and $\omega_-$.) 

Habiro has proved (again using quantum groups) that formulas of this type exist for all knots, but
the computation of the coefficients is not easy in general. Formulas of this kind are important for at least two reasons: computing quantities related to the Kashaev-Murakami-Murakami volume conjecture \cite{Ka,MM}, and computing Habiro's invariant of the homology sphere obtained by $\pm 1$ surgery on the knot. For more about this, see Habiro's survey article \cite{Hab2}.

{\bf Acknowledgements}\qua I got the idea for this work while
talking to T.Q.T.~Le during his visit to the University of Paris 7 in July 2002. I would 
like to thank both him and K.~Habiro for helpful discussions, and for sending me parts of their forthcoming papers \cite{Le2} and \cite{Hab3}.

\section{Habiro's formula for $\omega$}
 We use the notations of \cite{Hab} and of \cite{BHMV1}.  In particular, we write $$a=A^2, \ \ \{n\}=a^n-a^{-n}, \ \ [n]=\frac{a^n-a^{-n}}{a-a^{-1}}$$ and define $\{n\}!$ and $[n]!$ in the usual way.

Recall that the Kauffman bracket skein 
module, $\BK(M)$,  of an oriented $3$-manifold $M$ is the free $\BZ[A^{\pm}]$-module generated
by isotopy classes of banded links ($=$ disjointly embedded annuli) in $M$ modulo the submodule
generated by the Kauffman relations.  
\begin{figure}[ht!]
\begin{center}
$\KAUF \ = \ A \ \ \KAUFi \ + \ A^{-1}\ \ \KAUFii \ \ \ \ , \ \ \KAUFiii \ = -a-a^{-1}$
\caption{ The Kauffman bracket relations. (Recall $a=A^2$.)}\label{fig0}
\end{center}  
\end{figure}
  
 The Kauffman bracket gives an isomorphism $\langle \ \rangle:\BK(S^3)\mapright\approx \BZ[A^{\pm}]$. It is normalized so that the bracket of the empty link is $1$. 

The skein
module of the solid torus $S^1 \times D^2$ is $\BZ[A^{\pm}][z]$. We denote it by $\BB$. Here $z$ is  given by
the banded link $S^1\times J,$ where $J$ is
a small
arc in
 the interior of  
$D^2,$
and  $z^n$ means $n$ parallel copies of $z$. We define the even part $\BB^{ev}$ of $\BB$ to be the submodule generated by the even powers of $z$.

Let $t:\BB \rightarrow  \BB$  denote the twist map induced by a full right handed twist on the solid torus. It is well known 
(see {\em e.g.} \cite{BHMV1}) 
that there is a basis $ \{e_i\}_{i\ge0}$ of 
eigenvectors for the twist map. It is defined recursively by
 \begin{equation}\label{ei}
e_0=1, \quad e_1=z,  \quad e_i= z e_{i-1}- e_{i-2}~.
\end{equation}
 The eigenvalues are given by  
\begin{equation}\label{mui}t(e_i)=\mu_i e_i, \text { where } \mu_i=(-1)^i A^{i^2+2i} ~.\end{equation}

Let $\langle\ ,\ \rangle$ be the $\BZ[A^{\pm}]$-valued bilinear form on $\BB$ given by cabling the zero-framed Hopf link and taking the bracket.  For $x\in \BB$, put $\langle x\rangle=\langle x,1\rangle$. One has $\langle e_i\rangle=(-1)^i[i+1]$. 

Moreover, for every $f(z)\in \BB$, one has
\begin{equation}\label{lambdai} 
\langle f(z),e_i\rangle = f(\lambda_i)\langle e_i\rangle, \text{ where } \lambda_i= -a^{i+1}-a^{-i-1}~.
\end{equation}

Following Habiro \cite{Hab}, define $$R_n=\prod_{i=0}^{n-1} (z-\lambda_{2i}), \ \ S_n=\prod_{i=0}^{n-1} (z^2-\lambda_i^2)~.$$ The $R_n$ form a basis of $\BB$, and the $S_n$ form a basis of the even part $\BB^{ev}$ of $\BB$.

By construction, one has $\langle R_n, e_{2i}\rangle =0$ for $i<n$, and therefore also $\langle R_n, z^{2k}\rangle =0$ for $k<n$. Similarly, one has  $\langle S_n, e_{i}\rangle =0$ for $i<n$, and hence also $\langle S_n, z^{k}\rangle =0$ for $k<n$. It follows that $\langle R_n,S_m\rangle=0$ for $n\neq m$, and for $n=m$ one computes 
\begin{equation}\label{RS}
\langle R_n,S_n\rangle=\langle R_n,e_{2n}\rangle=\langle e_{2n}\rangle\prod_{i=0}^{n-1}(\lambda_{2n}-\lambda_{2i})=(-1)^n \frac{\{2n+1\}!}{\{1\}}~.
\end{equation}

We are looking for $$\omega_+=\sum_{n=0}^\infty c_{n,+} R_n$$ satisfying (\ref{fig-1}) for every even $x$, which is equivalent to requiring that 
\begin{equation}\label{zzz}
\langle \omega_+, x\rangle =\langle t(x)\rangle
\end{equation}
 for every $x\in \BB^{ev}$. Note that the left hand side of (\ref{zzz}) is a finite sum for every $x\in \BB^{ev}$.

\begin{thm}[Habiro\cite{Hab}] \label{0.1} Eq. (\ref{zzz}) holds  for 
\begin{equation}\label{Habf2}
c_{n,+}=(-1)^n \frac {a^{n(n+3)/2}}{\{n\}!}~.
\end{equation}
\end{thm}

Let us define $\omega_-$  to be the conjugate of $\omega_+$, where conjugation is defined, as usual, by $\overline A=A^{-1}$ and $\overline z=z$. Since conjugation corresponds to taking mirror images, we have that 
\begin{equation}\label{omega}\langle \omega_-, x\rangle =\langle t^{-1}(x)\rangle
\end{equation} for every even $x$.

 Note that $\omega_-=\sum_{n=0}^\infty c_{n,-} R_n$,  where 
\begin{equation}\label{Habf} c_{n,-}=\frac {a^{-n(n+3)/2}}{\{n\}!}~.
\end{equation}
 This follows from (\ref{Habf2}) since $\overline{R_n}=R_n$ and $\overline{ \{n\}}=-\{n\}$. 

\begin{rem}{\em The skein element $\omega=\omega_{+}$ is related to, but different from,  the element often called  $\omega$ appearing in the surgery axiom of Topological Quantum Field Theory (see for example \cite{BHMV2}).  If we call the latter $\omega^{TQFT}$, then  Equations (\ref{zzz}) and (\ref{omega}) would be  satisfied  by appropriate scalar multiples of $t^{-1}(\omega^{TQFT})$ and $t(\omega^{TQFT})$, respectively; moreover, they  would now hold not just for even $x$, but for all $x$. This applies in particular to the  $\omega$ of \cite{BHMV1}, which was constructed in a similar way as Habiro's $\omega_-$ (but using polynomials $Q_n=\prod_{i=0}^{n-1} (z-\lambda_i)$ in place of the polynomials $R_n$). 
}\end{rem} 
\section{A proof using orthogonal polynomials}

Habiro's proof of Theorem~\ref{0.1}  uses 
the relationship with the quantum group $U_q\fsl_2$.  Here is  another proof, using the method of orthogonal polynomials as in \cite{BHMV1}. 

Testing with the $S_n$-basis, we see that (\ref{omega}) holds if and only if 
$$ c_{n,-} = \frac {\langle t^{-1} S_n\rangle} {\langle R_n,S_n\rangle}~.$$ Thus, it is clear that an $\omega_-$ satisfying (\ref{omega}) exists, and to compute its coefficients, we just need to compute $\langle t^{-1} S_n\rangle$.

As in \cite{BHMV1}, define another bilinear form $\langle\ ,\ \rangle_1$ by $$\langle x, y\rangle_1=\langle t(x), t(y)\rangle~.$$ Define polynomials $\BR_n$ and $\BS_n$ by $$t(\BR_n)=\mu_n R_n, \ \ \ t(\BS_n)=\mu_{2n} S_n~.$$ (The factors $\mu_n$ and $\mu_{2n}$ are included so that $\BR_n$ and $\BS_n$ are monic, {\em i.e.} have leading coefficient equal to one.) 

Again, the $\BR_n$ form a basis of $\BB$, and the $\BS_n$ form a basis of the even part $\BB^{ev}$ of $\BB$, since the twist map $t$ preserves  $\BB^{ev}$. We have $\langle \BR_n,\BS_m\rangle_1=0$ for $n\neq m$, and 
\begin{equation}\label{curly}
\langle \BR_n,\BS_n\rangle_1=\mu_{n}\mu_{2n}\langle R_n,S_n\rangle~.
\end{equation}

Note that $\langle t^{-1} S_n\rangle=\mu_{2n}^{-1}\langle \BS_n\rangle$. Thus, we just need to compute $\langle \BS_n\rangle$.

\begin{prop} The polynomials $\BS_n$ satisfy a four-term recursion formula
\begin{equation}\label{rec}\BS_{n+1}=(z^2-\alpha_n)\BS_n-\beta_{n-1}\BS_{n-1}-\gamma_{n-2}\BS_{n-2}\end{equation} for certain  $\alpha_n, \beta_{n-1},\gamma_{n-2}\in \BZ[A^{\pm}].$
\end{prop}
\begin{proof} Since $\BS_n$ is monic of degree $2n$, we have that $z^2 \BS_n-\BS_{n+1}$ is a linear combination of the $\BS_k$ with $k\leq n$. The coefficients can be computed by taking  the scalar product with $\BR_k$. So we just need to show that  $\langle z^2 \BS_n, \BR_{k}\rangle_1=0$ if $k<n-2$. 

The point is that multiplication by $z$ is a self-adjoint operator with respect to the bilinear form $\langle\ ,\ \rangle_1$. In other words, one has $$\langle zx,y \rangle_1=\langle x , zy\rangle_1$$ for all $x,y\in\BB$. (This is because $\langle x,y \rangle_1=\langle t(xy) \rangle$.) It follows that 
$$\langle z^2 \BS_n, \BR_{k}\rangle_1=\langle \BS_n, z^2\BR_{k}\rangle_1=0 \ \ \text{if}\ \ k<n-2,$$ since $\BR_k$ has degree $k$, and $\BS_{n}$ annihilates all polynomials of degree $<n$.
\end{proof}
 
Note that the coefficients in the recursion formula (\ref{rec}) are given by 
\begin{equation}\label{all}
\alpha_n=\frac{\langle z^2 \BS_n, \BR_{n}\rangle_1}{\langle \BS_n, \BR_{n}\rangle_1},\ \beta_{n-1}=\frac{\langle z^2 \BS_n, \BR_{n-1}\rangle_1}{\langle \BS_{n-1}, \BR_{n-1}\rangle_1},\ \gamma_{n-2}=\frac{\langle z^2 \BS_n, \BR_{n-2}\rangle_1}{\langle \BS_{n-2}, \BR_{n-2}\rangle_1}~.
\end{equation} 

By convention, if $n<0$ then $\BR_n, \BS_n, \alpha_n, \beta_n,\gamma_n$ are all zero.

\begin{prop} One has
\begin{eqnarray}
\label{alpha} \alpha_n&=&2+a^{6n+4}[3]-a^{2n}\\
\label{beta} \beta_{n-1}&=&(a^{4n+1}+a^{8n+1}[3])\{2n\}\{2n+1\}\\
\label{gamma} \gamma_{n-2}&=&a^{10n-4}\{2n-2\}\{2n-1\}\{2n\}\{2n+1\}\end{eqnarray}
\end{prop}

\begin{proof} The formula for $\gamma_{n-2}$ is the easiest.  Let us use the notation $\bo_{\leq n}$ for terms of degree $\leq n$. Since $z^2\BR_{n-2}=\BR_n + \bo_{\leq n-1}$, we have $$\langle z^2 \BS_n, \BR_{n-2}\rangle_1 = \langle  \BS_n, z^2\BR_{n-2}\rangle_1 =\langle \BS_n, \BR_{n}\rangle_1~, $$ and hence formula  (\ref{gamma}) follows from (\ref{all}), (\ref{curly}), and (\ref{RS}).

For $\beta_{n-1}$, we need to compute
\begin{equation}\label{beta2}\langle z^2 \BS_n, \BR_{n-1}\rangle_1 = \langle  \BS_n, z^2\BR_{n-1}\rangle_1=\mu_{2n}\mu_{n-1} \langle  S_n, tz^2t^{-1}R_{n-1}\rangle~.
\end{equation}
This amounts to computing the coefficient of $R_n$ in the expression of $tz^2t^{-1}R_{n-1}$ in the $R_k$-basis. This coefficient can be computed as follows.

For $n\geq 1$, one has $$z^n=e_n +(n-1)e_{n-2} +\bo_{\leq n-4}~.$$ (This follows by induction from (\ref{ei}).) Thus, for $\varepsilon =\pm 1$, one has
$$t^\varepsilon z^n = \mu_n^\varepsilon z^n + (n-1) (\mu_{n-2}^\varepsilon - \mu_n^\varepsilon) z^{n-2} + \bo_{\leq n-4}~.$$ It follows that
\begin{equation}\label{tzt}
tz^2t^{-1} z^n = \frac {\mu_{n+2}} {\mu_{n}} z^{n+2} + (2-(n+1)\frac {\mu_{n+2}} {\mu_{n}} + (n-1) \frac {\mu_{n}} {\mu_{n-2}}) z^n +  \bo_{\leq n-2}~.
\end{equation}
Now write $$R_n=\prod_{i=0}^{n-1}(z-\lambda_{2i}) = z^n-x_{n-1}z^{n-1}+\bo_{\leq n-2}~, $$ where $x_{n-1}=\sum_{i=0}^{n-1} \lambda_{2i}$. Then (\ref{tzt}) gives
\begin{equation}\label{tztR}t z^2 t^{-1} R_{n-1} = \frac {\mu_{n+1}} {\mu_{n-1}} R_{n+1} + (x_n \frac {\mu_{n+1}} {\mu_{n-1}} - x_{n-2}\frac {\mu_{n}} {\mu_{n-2}}) R_n + \bo_{\leq n-1}
\end{equation}
and hence 
$$
\langle S_n, t z^2 t^{-1} R_{n-1} \rangle = (x_n \frac {\mu_{n+1}} {\mu_{n-1}} - x_{n-2}\frac {\mu_{n}} {\mu_{n-2}})\langle S_n , R_n\rangle~.$$
Plugging this into (\ref{beta2}), we have
\begin{align*}\langle z^2 \BS_n, \BR_{n-1}\rangle_1 =&\mu_{2n}\mu_{n-1} \langle  S_n, tz^2t^{-1}R_{n-1}\rangle\\
=&(x_n \frac {\mu_{n+1}} {\mu_{n}} - x_{n-2}\frac {\mu_{n-1}} {\mu_{n-2}})\langle \BS_n , \BR_n\rangle_1\\
=& (A^{6n+1} [3] +A^{-2n+1})\langle \BS_n , \BR_n\rangle_1~.
\end{align*}
Using  (\ref{all}), (\ref{curly}), and (\ref{RS}) as before, this implies formula (\ref{beta}) for $\beta_{n-1}$.

Finally, for $\alpha_n$, let us compute
\begin{equation}\label{alpha2}\langle z^2 \BS_n, \BR_{n}\rangle_1 = \mu_{2n}\mu_{n} \langle  tz^2t^{-1} S_n, R_{n}\rangle~.
\end{equation}
This amounts to computing the coefficient of $S_n$ in the expression of $tz^2t^{-1}S_{n}$ in the $S_k$-basis.\footnote{This is easier than computing $\langle   S_n, tz^2t^{-1}R_{n}\rangle$ by expanding $tz^2t^{-1}R_{n}$ in the $R_k$-basis, because the latter would require computing the first {\em three} terms, and not just the first two terms as in (\ref{tztR}) above and also in (\ref{tztS}) below.} 

The computation is similar to the one above. We write $$S_n=\prod_{i=0}^{n-1}(
z^2-\lambda_{i}^2) = z^{2n}-y_{n-1}z^{2n-2}+\bo_{\leq 2n-4}~, $$ where $y_{n-1}=\sum_{i=0}^{n-1} \lambda_{i}^2$. Then (\ref{tzt}) gives
\begin{align}\label{tztS}t z^2 t^{-1} &S_{n} =\\ \frac {\mu_{2n+2}} {\mu_{2n}} S_{n+1} + &\left(2+ (y_n-2n-1) \frac {\mu_{2n+2}} {\mu_{2n}} - (y_{n-1}-2n+1)\frac {\mu_{2n}} {\mu_{2n-2}}\right) S_n + \bo_{\leq 2n-2}\nonumber
\end{align}and hence 
\begin{align*}
\langle t z^2 t^{-1}S_n, R_{n} \rangle =& \left(2+ (y_n-2n-1) \frac {\mu_{2n+2}} {\mu_{2n}} - (y_{n-1}-2n+1)\frac {\mu_{2n}} {\mu_{2n-2}}\right)\langle S_n , R_n\rangle\\
=&\,\, (2+a^{6n+4}[3]-a^{2n})\langle S_n , R_n\rangle~. 
\end{align*}
Plugging this into (\ref{alpha2}), we get 
$$\langle z^2 \BS_n, \BR_{n}\rangle_1 = (2+a^{6n+4}[3]-a^{2n}) \langle  \BS_n, \BR_{n}\rangle_1~,$$ proving formula (\ref{alpha}) for $\alpha_{n}$.
\end{proof}

\begin{proof}[Proof of Habiro's Theorem \ref{0.1}] As already observed, we have $$c_{n,-}=\frac {\langle t^{-1} S_n\rangle} {\langle R_n,S_n\rangle}= \mu_{2n}^{-1} \frac {\langle \BS_n\rangle} {\langle R_n,S_n\rangle}~.$$
But $\langle \BS_n\rangle$ satisfies the recursion relation
\begin{equation}\label{rec2}\langle \BS_{n+1}\rangle =(\lambda_0^2-\alpha_n)\langle \BS_n\rangle -\beta_{n-1}\langle \BS_{n-1}\rangle -\gamma_{n-2}\langle \BS_{n-2}\rangle
\end{equation} (since $\langle z\rangle=\lambda_0$). It follows that 
\begin{equation}\label{rec3} \langle \BS_n\rangle = (-1)^n a^{(3n^2+n)/2} \frac {\{n+1\}\{n+2\}\cdots\{2n+1\}}{\{1\}}~, \end{equation} since one can check that (\ref{rec3}) is true for $n=0,1,2$ and that it solves the recursion (\ref{rec2}). This  implies Habiro's formula (\ref{Habf}) for $c_{n,-}$. Taking conjugates, one then also obtains formula (\ref{Habf2}) for $c_{n,+}$.
\end{proof}
\begin{rem}{\em Although it might be hard to guess formula (\ref{rec3}), once one knows it the recursion relation (\ref{rec2}) is easily checked.  Observe that $\lambda_0^2=a^2+a^{-2}-2$. Put $q(n)=(3n^2+n)/2$. Then (\ref{rec2}) is equivalent to  
\begin{align*} (a^2+a^{-2} &-a^{6n+4}[3]+a^{2n}) a^{q(n)} +(a^{4n+1}+a^{8n+1}[3]) \{n\} a^{q(n-1)}\\&- a^{10n-4} \{n-1\}\{n\} a^{q(n-2)} 
= -  \frac {\{2n+2\}\{2n+3\}} {\{n+1\}}a^{q(n+1)}
\end{align*} which is a straightforward computation.
} 
\end{rem}

\section{Graphical calculus and a formula for $\omega^p$}\label{sec3}

Let us write $\omega=\omega_+$ and put 
\begin{equation}\label{om}
\omega^p= \sum_{n=0}^\infty c_{n,p} R_n~.
\end{equation}
 Note that the coefficients $c_{n,p}$ are well-defined (because  $R_n$ divides $R_{n+1}$ and therefore the coefficients $C_{n,m}^k$ in the product expansion $R_n R_m = \sum_{k} C_{n,m}^k R_k$ are zero if $n$ or $m$ is bigger than $k$.) We have
\begin{equation}\label{omp}
\langle \omega^p, x\rangle =\langle t^p(x)\rangle
\end{equation}
 for every even $x$. (This follows from (\ref{zzz}) since circling with $\omega^p$ is the same as circling with $p$ parallel copies of $\omega$.) Of course,  $c_{n,1}=c_{n,+} $ and $c_{n,-1}=c_{n,-} $ (it follows from the uniqueness of $\omega_-$ that $\omega_-= \omega^{-1}$). The aim of this section is to give a formula for the coefficients $c_{n,p}$ (see Theorem~\ref{th2} below).

\begin{figure}[ht!]
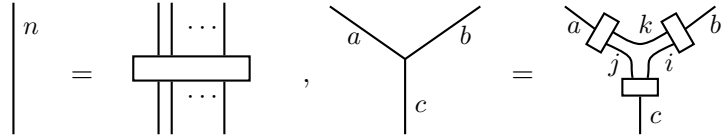

\begin{center}
\LINE \ \ $=$\ \ \BOX\ \ \ , \VERTEX \ \ $=$\ \ \VERTEXi
\caption{\label{fig1} How to expand colored edges and vertices. The boxes stand for appropriate Jones-Wenzl idempotents.}  
\end{center}
\end{figure}

We use the extension of the Kauffman bracket to admissibly colored banded trivalent graphs as in \cite{MV}. (Such graphs are sometimes called spin networks; for more background see {\em e.g.} \cite{KL} and references therein.) A color is just an integer $\geq 0$. A triple of colors $(a,b,c)$ is {\em admissible} if $a+b+c \equiv 0 \pmod 2$ and $|a-b|\leq c\leq a+b$. Let $D$ be a planar diagram of a banded trivalent graph.  An admissible coloring of $D$ is an assignment of colors to the edges of $D$ so that at each vertex, the three colors meeting there form an admissible triple. The Kauffman bracket of $D$ is defined to be the bracket of the {\em expansion} of $D$ obtained as follows. The expansion of an edge colored $n$ consists of  $n$ parallel strands with a copy of the Jones-Wenzl idempotent $f_n$ inserted. (The idempotent $f_n$ is characterized by the fact
that $x f_n= f_n x=0$ for every element $x$ of the standard basis of the Temperley-Lieb algebra other than the identity element; here, the standard basis consists of the $(n,n)$-tangle diagrams without crossings 
and without closed loops.) The expansion of a vertex is defined as in Fig. \ref{fig1}, where the {\em internal colors} $i,j,k$ are defined by
\begin{equation} \label{intco} i=(b+c-a)/2, \ j=(c+a-b)/2, \ k=(a+b-c)/2~.
\end{equation}

We have the {\em fusion} equation
\begin{equation}\label{fus}
\FUS \ \ \   = \sum_c \frac {\langle c \rangle}{\langle a,b,c\rangle} \  \ \ \  \FUSi
\vspace{.2cm}
\end{equation} Here the sum is over those colors $c$ so that the triple $(a,b,c)$ is admissible, we have $\langle c \rangle=\langle e_c \rangle=(-1)^c[c+1]$,  and the trihedron coefficient $\langle a,b,c\rangle$ is (see \cite[Thm.~1]{MV}):
\begin{equation}\label{theta}
\langle a,b,c\rangle = \THETA =(-1)^{i+j+k}
\frac {[i+j+k+1]!\,[i]!\,[j]!\,[k]!}{[a]!\,[b]!\,[c]!}
\vspace{.2cm}
\end{equation}
(here $i,j,k$ are the internal colors as defined in  (\ref{intco})).
Note that $\langle n,n,2n\rangle=\langle 2n\rangle $ so that
\begin{equation}\label{fusa}
\FUSa \ \ \   =  \  \ \ \  \FUSai\ \ \ \  + \ \ \ldots
\vspace{.4cm}
\end{equation}

We will need the following lemma.
\begin{lem} \label{lem1} For $0\leq k\leq n$,  one has
$$\CAL  = \frac {([k]!)^2}{[2k]!} \CALo$$
\end{lem}
\begin{proof} This follows from the formula for the tetrahedron coefficient given in \cite[Thm. 2]{MV}. (The  sum over $\zeta$ in that formula reduces to just one term.) \end{proof}

The key observation is that 
\begin{equation} \label{eq3}
\CALiia \ \   = \ 0 \ \ \ \text{for} \ \ k\neq n~.
\end{equation}
Indeed, if $k<n$, there are at most $2k$ vertical strands in the middle and so the result is zero because of the Jones-Wenzl idempotents $f_{2n}$ at the top and bottom.  On the other hand, if $k>n$, the result is zero because $R_k$ annihilates all even polynomials in $z$ of degree $<2k$. 

If $k=n$ one finds
\begin{equation} \label{eq4}
\CALii = (-1)^n (\{n\}!)^2 \CALo~.
\end{equation}

Indeed, applying (\ref{fusa}), one has 
$$\CALii = \CALiii = \ \frac {\langle R_n,e_{2n}\rangle }{\langle e_{2n}\rangle}
\ \ \ \CALa $$
hence (\ref{eq4}) follows from (\ref{RS}) and Lemma \ref{lem1}.

On the other hand, since circling with $\omega^p$ induces $p$ full twists on even numbers of strands (see (\ref{omp})), we have, using  (\ref{eq3}), that

\begin{equation} \label{eq5}\mu_n^{-2p} c_{n,p} \CALii \ = \mu_n^{-2p}  \CALiic = \ \CALiv  
\end{equation}
where there are $2p$ crossings  in the last diagram. (See (\ref{mui}) for the twist eigenvalues $\mu_i$.) Applying the fusion equation (\ref{fus}), we have
\begin{equation}\label{improve} \CALiv = \sum_{k=0}^n \delta(2k; n,n)^{2p} \frac {\langle 2k \rangle}{\langle n,n,2k\rangle}  \ \ \ \CAL
\end{equation} where $\delta(c;a,b)$ is the {\em half-twist coefficient} defined by $$\DEL \ \ \   = \delta(c;a,b) \  \ \ \  \DELi~.$$ This coefficient is computed in \cite[Thm. 3]{MV}. For us, it is enough to know that $$\delta(c;a,b)^2=\frac{\mu_c}{\mu_a\mu_b}$$ which is easy to see. Using (\ref{eq4}) and Lemma \ref{lem1}, it follows that 

$$\mu_n^{-2p} c_{n,p} (-1)^n (\{n\}!)^2 = \sum_{k=0}^n \frac{\mu_{2k}^p}{\mu_n^{2p}} \frac {\langle 2k \rangle}{\langle n,n,2k\rangle} \frac {([k]!)^2}{[2k]!}~. $$ The factors of $\mu_n^{-2p}$ cancel out, and in view of (\ref{theta}), this gives the following result:
\begin{thm}\label{th2} The coefficients $c_{n,p}$ of $\omega^p$ in (\ref{om}) are given by 
\begin{equation} \label{cnp} c_{n,p} = \frac{1}{(a-a^{-1})^{2n}} \sum_{k=0}^n  \frac {(-1)^k \mu_{2k}^p [2k+1]}{[n+k+1]!\, [n-k]!}~.
\end{equation}
\end{thm}

\section{Another formula for the coefficients of $\omega^p$}
 
Following Habiro, we introduce the polynomials $R_n'= (\{n\}!)^{-1} R_n$ and write 
\begin{equation}\label{omprim}
\omega^p= \sum_{n=0}^\infty c'_{n,p} R_n'~,
\end{equation}
 where $$c'_{n,p}= \{n\}! \,c_{n,p}~.$$ The aim of this section is to show that $c'_{n,p}$ is a Laurent polynomial, {\em i.e.} that $c'_{n,p}\in \BZ[A^{\pm}]$. This fact was shown by Habiro \cite{Hab3} using the quantum group $U_q\fsl_2$. 
Observe that by (\ref{Habf}) and (\ref{Habf2}), we already know this 
fact for $p=\pm 1$:

$$  c'_{n,1} =(-1)^n a^{n(n+3)/2}~, \ \  c'_{n,-1} = a^{-n(n+3)/2}~.$$

(But  Formula (\ref{cnp})  tells us only that 
\begin{equation} \label{cnp'} c'_{n,p} = \frac{1}{(a-a^{-1})^{n}} \sum_{k=0}^n (-1)^k \mu_{2k}^p [2k+1] \frac {[n]!}{[n+k+1]!\, [n-k]!}~,
\end{equation} from which it is not clear that $c'_{n,p}\in \BZ[A^{\pm}]$.)

To do so, we will replace Formula (\ref{improve}) in the previous section by Formula (\ref{xxx}) below. For this, we need the
following two Lemmas. 

\begin{lem} We have 
$$\REC \ \ = \ \sum_{k=0}^n C_{n,k} \RECa~, $$ where 
\begin{equation}\label{cnk} C_{n,k}= a^{n(n-k)} \qchoose{n} {k} \prod_{j=n-k+1}^n (1-a^{-2j})~.
\end{equation}
\end{lem} Here, as usual, $$\qchoose{n} {k} = \frac {[n] [n-1] \cdots [n-k+1]}{[k]!}~.$$

\begin{proof} For $0\leq p\leq n$, let us write, more generally,
  $$\RECj \ \ = \sum_{k=0}^n \, C_{n,p,k} \RECaj~. $$ First, we  consider the case $p=1$. 
By induction on $n$, it is easy to prove that

$$\RECi \ \ = \ a^n \RECii + a^{n-1}(1-a^{-2n}) \RECiii$$
(recall $a=A^2$). Using this, we now fix $n$ and do induction on $p$ to 
obtain 
the recursion formula 
\begin{equation}\label{recf} C_{n,p+1,k}= a^{n-2k}C_{n,p,k} + a^{n-2k+1} (1-a^{-2(n-k+1)}) C_{n,p,k-1}~.
\end{equation} 
Here we have used the following two facts which follow from the defining properties of the Jones-Wenzl idempotents.\footnote{Equation (\ref{Id2}) is a special case of the half-twist coefficient.} 

\begin{align}\label{Id1} \IDEM \ \ \ &=\  \ \ \IDEMiii \\
\label{Id2}\IDEMii \ \ \ &=\ \ A^{-pq}\ \ \ \ \IDEMiii
\end{align}  

Note that the coefficients $C_{n,p,k}$ behave like the binomial 
coefficients $\mychoose{p}{k}$ in that $C_{n,0,0}=1$, and $C_{n,p,k}=0$ for $k<0$ or $k>p$. It follows that the recursion formula (\ref{recf}) determines the $C_{n,p,k}$ uniquely. One finds
$$C_{n,p,k}= a^{p(n-k)} \qchoose{p} {k} \prod_{j=n-k+1}^n (1-a^{-2j})~.$$ Specializing to the case $p=n$, this proves the Lemma.
\end{proof}

\begin{rem}{\em The coefficient $C_{n,n}$ was computed by a different method in \cite[Prop. 4.4]{Abchir}. Knowing this coefficient would be enough to obtain Habiro's formula (\ref{Habf2}) for $\omega$ (see Remark~\ref{p=1} below). Unfortunately, the method of \cite{Abchir} does not give the coefficients $C_{n,k}$ for $k\neq n$, which we will need to obtain a formula for $\omega^p$.
}\end{rem}

\begin{lem} We have
\begin{equation}\label{diaglem} \RECb \ \ = \frac{\mu_{n-k}^2}{\mu_n^2} \RECc 
\end{equation}
\end{lem}

\begin{proof} The left hand side of (\ref{diaglem}) is equal to 

$$\RECd $$ By an isotopy, this becomes $$\mu_k^{-2} \RECe $$

Applying the coefficients $\delta(n;n-k,k)^{-2}$ and $\delta(n;k,n-k)^{-2}$, which are both equal to $\mu_k\mu_{n-k}\mu_n^{-1}$ (see also (\ref{Id2})), we see that this is equal to the right hand side of (\ref{diaglem}).
\end{proof}

Let $C_{n,k}^{(p)}$ be the coefficient defined by the expansion
\begin{equation}\label{pcross}
\RECp \ \ =\sum_{k=0}^n \,C_{n,k}^{(p)} \RECa~,
\end{equation} (where the diagram on the left hand side of (\ref{pcross}) has $2p$ crossings).  Putting the two preceding Lemmas together, we may obtain a formula for this coefficient.  In particular, it follows by induction on $p$ that $C_{n,k}^{(p)}$ is a Laurent polynomial divisible by $\prod_{j=n-k+1}^n (1-a^{2j})$. 

We are interested in the coefficient  $C_{n,n}^{(p)}$, since we have
\begin{equation}\label{xxx}
\CALiv \ = \ C_{n,n}^{(p)}\ \ \  \CALip \ = \ C_{n,n}^{(p)} \CALo
\end{equation}
(where there are $2p$ crossings in the diagram on the left). 
Using (\ref{eq4}) and (\ref{eq5}) from Section~\ref{sec3}, it follows  that 

$$ \mu_n^{-2p} c_{n,p} (-1)^n (\{n\}!)^2 = C_{n,n}^{(p)} $$
and therefore 
\begin{equation} \label{cprim} c'_{n,p}= \{n\}!\, c_{n,p} = (-1)^n \mu_n^{2p}(\{n\}!)^{-1}C_{n,n}^{(p)}~.\end{equation} As already observed,  $C_{n,n}^{(p)}$ is divisible by $$\prod_{j=1}^n (1-a^{2j})= a^{-n(n+1)/2} \{n\}!$$ and hence $c'_{n,p}$ is indeed a Laurent polynomial.

\begin{rem}\label{p=1} {\em In the case $p=1$, we have 
\begin{equation}\label{yyy}
C_{n,n}^{(1)}= C_{n,n} = a^{-n(n+1)/2} \{n\}!~.
\end{equation}
Plugging this into (\ref{cprim}), we get $$c'_{n,1}=(-1)^n a^{n^2+2n}a^{-n(n+1)/2}= (-1)^n a^{n(n+3)/2}~,$$   giving another proof of Habiro's formula (\ref{Habf2}) for $c_{n,1}$. 
}\end{rem}

Here is an explicit formula for $C_{n,n}^{(p)}$ which follows from (\ref{cnk}) and (\ref{diaglem}). The sum is over all multi-indices $\underline k =(k_1,\ldots, k_p)$ such that $k_i\geq 0$ for all $i$, and $\sum k_i = n$. For convenience, put $s_i=k_1+\ldots + k_i$  and $r_i=n-s_i=k_{i+1}+\ldots+k_p$, and define
\begin{equation}\label{phik}\varphi(\underline{k})= \sum_{i=1}^{p-1} 
r_i(r_{i-1}+r_i+2)~.
\end{equation}

\begin{align*} C_{n,n}^{(p)}&= \sum_{\underline k} C_{n,k_1} \frac{\mu_{n-k_1}^2}{\mu_n^2} C_{n-k_1, k_2}\frac{\mu_{n-k_1-k_2}^2}{\mu_{n}^2} \cdots 
C_{n-s_{p-2}, k_{p-1}}\frac{\mu_{n-s_{p-1}}^2}{\mu_{n}^2}C_{n-s_{p-1}, k_{p}}\\ 
&=\mu_n^{2-2p}\sum_{\underline k} \left(\prod_{i=1}^{p-1}\mu_{r_i}^2\right) \prod_{i=1}^{p} C_{r_{i-1}, k_i}\\
&= \mu_n^{2-2p} \left(\prod_{j=1}^n (1-a^{-2j})\right)\sum_{\underline k} \left(\prod_{i=1}^{p-1}a^{{r_i}^2+2r_i}\right) \prod_{i=1}^{p-1} \left(a^{r_{i-1}r_i} \qchoose{r_{i-1}}{k_i}\right)\\
&= \mu_n^{2-2p} a^{-n(n+1)/2} \{n\}! \sum_{\underline k}  a^{\varphi(\underline{k})} \qchoose{n}{\,\,\underline{k}\,\,}
\end{align*} 
where we have put, as usual, $${\qchoose{n}{\,\,\underline{k}\,\,}}=\frac{[n]!}{[k_1]! \cdots [k_p]!}~.$$

In view of (\ref{cprim}), and using $\mu_n^2=a^{n^2+2n}$, we obtain the following final result:
\begin{thm}[Habiro \cite{Hab3}]\label{th3} For $p\geq 1$, the coefficients $c'_{n,p}$ of $\omega^p$ in (\ref{omprim}) are given by  
\begin{equation} c'_{n,p}= (-1)^n a^{n(n+3)/2}\sum_{\underline k =(k_1,\ldots,k_p) \atop k_i\geq 0, \ \sum k_i = n}  a^{\varphi(\underline{k})} \qchoose{n}{\,\,\underline{k}\,\,}~.\end{equation}
where  $\varphi(\underline k)$ is defined in (\ref{phik}). 
\end{thm}

\begin{rem}{\em Since $c'_{n,-p}=(-1)^n \overline{c'_{n,p}}$, this also determines the coefficients of negative powers of $\omega$.}\end{rem}

\begin{rem}{\em In \cite{Hab3}, Habiro has obtained a similar formula using 
the quantum group $U_q\fsl_2$. 
}\end{rem}
\noindent{\bf Example:}  Assume  $p=2$. We may write $\underline{k}=(k,n-k)$. Then \begin{align*}  c'_{n,2}&= (-1)^n a^{n(n+3)/2}\sum_{k=0}^n a^{(n-k)(n+n-k +2)}  \qchoose nk\\
&=(-1)^n a^{(5n^2+7n)/2} \sum_{k=0}^n a^{k^2-2k-3nk} \qchoose nk
\end{align*}

\section{The colored Jones polynomial of twist knots}\label{twistknots}

\begin{figure}[ht!]
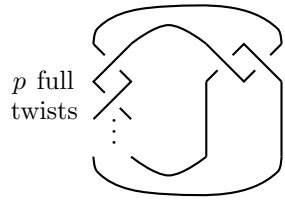

\begin{center}
\TW 
\caption{ The twist knot $K_p$. (Here $p\in\BZ$.) For $p=1$, $K_p$ is a left-handed trefoil, and for $p=-1$, $K_p$ is the figure eight knot.} \label{fig4} 
\end{center}
\end{figure}

In this last section, we illustrate one use of $\omega^p$, namely to give a formula for the colored Jones polynomial of twist knots (see Figure~\ref{fig4}) in terms of the coefficients $c'_{n,p}$. The results of this section are known to K.~Habiro and T.Q.T.~Le.

The colored Jones polynomial of a knot $K$ colored with the $N$-dimensional irreducible representation of $\fsl_2$ can be expressed as the Kauffman bracket of $K$ cabled by $(-1)^{N-1}e_{N-1}$: $$J_K(N)= (-1)^{N-1} \langle K(e_{N-1})\rangle ~.$$ (The factor of $(-1)^{N-1}$ is included so that $J_{Unknot}(N)=[N]$.) We will use the normalization  $$J'_K(N)=\frac {J_K(N)}{J_{Unknot}(N)}=\frac{ \langle K(e_{N-1})\rangle}{\langle e_{N-1}\rangle}~.$$ Here, we assume the knot $K$ is equipped with the zero framing.

\begin{figure}[ht!]
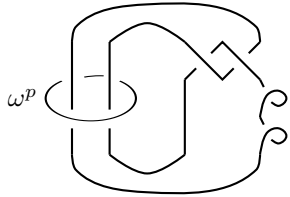

\begin{center}
$$\TWa$$
\caption{A surgery description of the twist knot $K_p$ with zero framing.} \label{fig5} 
\end{center}
\end{figure}

Let us compute $\langle K_p(e_{N-1})\rangle$. We use the surgery description given in Fig.~\ref{fig5}. Recall that $\omega^p=\sum c'_{k,p} R'_k$. By induction, one can check that   
\begin{equation}\label{eR}
e_{N-1}=\sum_{n=0}^{N-1} (-1)^{N-1-n} \qchoose{N+n}{N-1-n} R_n~.
\end{equation}
The key observation (which I learned from T.Q.T.~Le) is that $$\TWan \ \ \ = \ \ 0 $$ for $k\neq n$. This is because each component of this link is a zero-framed unknot having a spanning disk pierced twice by the other component, and circling with $R_m$ annihilates all even polynomials in $z$ of degree $<2m$.

Thus only terms with $k=n$ survive, and so we have
\begin{equation}\label{coli}
\langle K_p(e_{N-1})\rangle \ = \ \sum_{n=0}^{N-1} (-1)^{N-1-n} \qchoose{N+n}{N-1-n}  c'_{n,p} \TWnn
\end{equation}
\vspace{8pt}
Now, using that $R_n-e_n$ has degree $<n$, we compute
\begin{align*} &\TWnn \ \ \ = \TWnnn  \\ 
&= \ \mu_n^2 \TWi \ \ \ = \ \mu_n^2  C_{n,n} \langle R'_n,e_{2n}\rangle \\
&= (-1)^n a^{n(n+3)/2}\frac{\{2n+1\}!}{\{1\}} 
\end{align*}
where we have used (\ref{xxx}) in the last but one equation, and $C_{n,n}$ and ${\langle R_n,e_{2n}\rangle }$ are given in (\ref{yyy}) and (\ref{RS}), respectively (but notice we are using $R_n'=(\{n\}!)^{-1}R_n$ here).

Plugging this into (\ref{coli}), we obtain the following result.
\begin{thm} The colored Jones polynomial of the twist knot $K_p$ is given by \begin{equation} \label{col} J'_{K_p}(N)= \sum_{n=0}^{\infty}  f_{K_p,n} \frac{\{N-n\} \{N-n+1\} \cdots \{N+n\}}{ \{N\}} ~,
\end{equation} where $$f_{K_p,n} =  a^{n(n+3)/2} c'_{n,p}~.$$ (The sum is actually finite,  the terms with $n\geq N$ being all zero.)
\end{thm}
 Note that since $c'_{n,p}$ is a Laurent polynomial in $a$, so are the coefficients $f_{K_p,n}$. For example, for the figure eight knot $K_{-1}$, we have $f_{K_{-1},n}=1$, and for the left-handed trefoil $K_1$, we have $f_{K_{1},n}=(-1)^n  a^{n(n+3)}$. (The right-handed trefoil is the mirror image of $K_1$, so one just needs to take the conjugate of  $J'_{K_1}(N)$.)  These formulas can be found in \cite{Hab2}. 
\begin{rem}{\em Here is another expression for $J'_{K_p}(N)$. Put $q=a^2=A^4$ and, as is customary in $q$-calculus, $$(x)_n=(1-x)(1-xq)\cdots (1-xq^{n-1})~.$$ Then (\ref{col}) gives 
$$J'_{K_p}(N)= \sum_{n=0}^{\infty}  \tilde f_{K_p,n} (q^{1-N})_n (q^{1+N})_n ~,$$ where $$\tilde f_{K_p,n}= (-1)^n q^{-n(n+1)/2} f_{K_p,n}~.$$ For example, for the figure eight knot $K_{-1}$, we have $\tilde f_{K_{-1},n}=(-1)^n q^{-n(n+1)/2}$, and for the left-handed trefoil $K_1$, we have $\tilde f_{K_{1},n}=q^n$. These formulas can be found in \cite{Le2} (see also \cite{Le1}).
}\end{rem}

\Addressesr


\begin{thebibliography}{BHMV2}

\bibitem[A]{Abchir}{\sc H. Abchir.} TQFT invariants at infinity for the Whitehead manifold. In: {\em Knots in Hellas '98}, Proceedings of the International
     Conference on Knot Theory and Its
                               Ramifications, World Scientific, Series on Knots and Everything, Vol. 24. 
\bibitem[BHMV1]{BHMV1} {\sc  C. ~Blanchet, N. ~Habegger, G. ~Masbaum,  P. ~Vogel.}
Three-manifold invariants derived from the Kauffman bracket,  {\em Topology} {\bf 31} (1992), 685-699
\bibitem[BHMV2]{BHMV2} {\sc  C. ~Blanchet, N. ~Habegger, G. ~Masbaum,  P. ~Vogel.} Topological quantum field theories derived from the Kauffman bracket,  {\em Topology}   {\bf 34} (1995), 883-927

\bibitem[H1]{Hab} {\sc K. Habiro.} On the colored Jones polynomial of
  some simple links. In: Recent Progress Towards the Volume
  Conjecture, Research Institute for Mathematical Sciences (RIMS) Kokyuroku 1172, September 2000.

\bibitem[H2]{Hab2}  {\sc K. Habiro.}  On the quantum $sl_2$ invariants of knots and integral homology spheres. In: Invariants of knots and 3-manifolds
   (Kyoto 2001), Geometry and Topology Monographs, Vol. 4 (2002), 55-68. 
\bibitem[H3]{Hab3}{\sc K. Habiro.} In preparation.

\bibitem[Ka]{Ka}{\sc R.M. Kashaev.} The hyperbolic volume of knots from quantum dilogarithm. {\em Lett. Math. Phys.} {\bf 39} (1997) 269-275.

 \bibitem[KL]{KL} {\sc L. H. Kauffman, S. L. Lins.} {\em Temperley-Lieb recoupling theory and invariants of $3$-manifolds,} Ann.\ Math.\ Studies 133 (Princeton University Press, 1994). 

\bibitem[L1]{Le1} {\sc T.Q.T. Le.}  Quantum invariants of 3-manifolds: integrality, splitting, and perturbative expansion. {\tt arXiv:math.QA/0004099}

\bibitem[L2]{Le2} {\sc T.Q.T. Le.}  In preparation.


\bibitem[MV]{MV} {\sc G. Masbaum, P. Vogel.} $3$-valent graphs and the Kauffman
bracket.  {\em Pacific J.\ Math.} {\bf 164}, (1994) 361-381.
\bibitem[MM]{MM} {\sc H. Murakami, J. Murakami.} The colored Jones polynomial and the simplicial volume of a knot. {\em Acta Math.} {\bf 186} (2001) 85-104.
\end{thebibliography}
\end{document}